\documentclass[american,sort&compress, 12pt]{article}
\usepackage[T1]{fontenc}
\usepackage[latin9]{inputenc}
\usepackage[a4paper]{geometry}
\usepackage{units}
\usepackage{amsmath}
\usepackage{amssymb}
\usepackage{graphicx}
\usepackage[numbers]{natbib}

\makeatletter
\newcommand{\lyxaddress}[1]{
\par {\raggedright #1
\vspace{1.4em}
\noindent\par}
}

\usepackage{hyperref}
\hypersetup{
    colorlinks,
    citecolor=black,
    filecolor=black,
    linkcolor=black,
    urlcolor=black
}

\usepackage{appendix}
\usepackage{chngcntr}
\usepackage{etoolbox}

\AtBeginEnvironment{subappendices}{%
\section*{Appendix}
\addcontentsline{toc}{section}{Appendices}
}

\makeatother

\usepackage{babel}
\begin{document}

\title{The dynamical impact of a shortcut in unidirectionally coupled rings
of oscillators}

\author{Jan Philipp Pade$^{1}$, Leonhard Lücken$^{1}$ and Serhiy Yanchuk$^{1}$}
\date{}
\maketitle

\lyxaddress{$^{1}$Humboldt-University of Berlin, Institute of Mathematics, Unter
den Linden 6, 10099 Berlin, Germany}
\begin{abstract}
We study the destabilization mechanism in a unidirectional ring of
identical oscillators, perturbed by the introduction of a long-range
connection. It is known that for a homogeneous, unidirectional ring
of identical Stuart-Landau oscillators the trivial equilibrium undergoes
a sequence of Hopf bifurcations eventually leading to the coexistence
of multiple stable periodic states resembling the Eckhaus scenario.
We show that this destabilization scenario persists under small non-local
perturbations. In this case, the Eckhaus line is modulated according
to certain resonance conditions. In the case when the shortcut is
strong, we show that the coexisting periodic solutions split up into
two groups. The first group consists of orbits which are unstable
for all parameter values, while the other one shows the classical
Eckhaus behavior.
\end{abstract}

\section{Introduction}

The theory of dynamics on networks is a rapidly growing field of research.
Motivated by applications from neuroscience, biological systems, and
social networks, scientists from various areas have investigated the
connection between network structure and dynamics in the last few
decades. Since then, the importance of special topologies like scale-free
or small world networks was discovered and dynamics on these structures
were studied extensively \citep{Boccaletti2006}.

\begin{figure}
\begin{centering}
\includegraphics[width=0.3\textwidth]{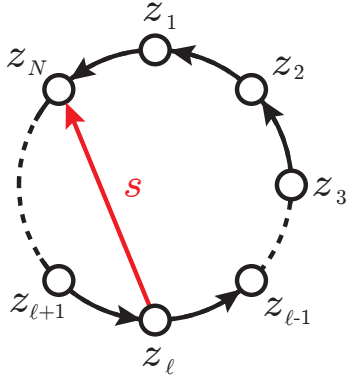}
\par\end{centering}

\caption{\emph{\label{fig:Coupling-scheme}}Coupling scheme of a unidirectional
ring with a shortcut.}
\end{figure}

In this article, we study a special class of network topologies, namely
unidirectional rings which are perturbed by the insertion of a single
additional non-local link {[}Fig.~\ref{fig:Coupling-scheme}{]}.
One reason for the interest in ring structures is that they emerge
in many natural systems \citep{Koseska2010,Restrepo2004a,Takamatsu2001}.
Also, rings can be seen as motifs of larger and more complex networks
\citep{Milo2002}. A lot of research has been done for bidirectional
rings in the last few years \citep{Abrams2004,Daido1997,Restrepo2004,Waller1984,Zou2009}.
In contrary, less is known about dynamics in unidirectional rings
\citep{horikawa2012b,Horikawa2012a,Perlikowski2010a,Yanchuk2008a,Popovych2011,Perlikowski2010},
although these structures play an important role in various applications
\citep{Bressloff1997,collins1994,strelkowa2011,Takamatsu2001,vishwanathan2011,VanderSande2008}.
As a simple, paradigmatic model, we consider unidirectional rings
of $N$ identical Stuart-Landau oscillators with an additional shortcut
from node $\ell$ to node $N$, to which we refer as a perturbation
of the homogeneous system. The dynamics on the perturbed ring are
described by the following equations: 
\begin{eqnarray}
\dot{z}_{j}\left(t\right) & = & \left(\mu-\left|z_{j}\left(t\right)\right|^{2}\right)z_{j}+z_{j+1}\left(t\right),\quad j=1,...,N-1,\nonumber \\
\dot{z}_{N}\left(t\right) & = & \left(\mu-\left|z_{N}\left(t\right)\right|^{2}\right)z_{N}+z_{1}\left(t\right)+sz_{\ell}\left(t\right),\label{eq:SL-system}
\end{eqnarray}
where $\mu=\alpha+i\beta$, $\beta>0$, $s>0$ is the shortcut strength,
and $z_{j}\left(t\right)\in\mathbb{C}$. For an unperturbed ring ($s=0$),
the destabilization mechanism is described in \citep{Yanchuk2008a}.
In particular, when the parameter $\alpha=\Re\left(\mu\right)$ in
(\ref{eq:SL-system}) is increased, the zero equilibrium 
\begin{equation}
z_{1}=...=z_{N}=0\label{eq:zero-equilibrium}
\end{equation}
loses stability and undergoes a sequence of Hopf bifurcations. The
first half of the emerging branches of periodic solutions stabilizes
at appropriate values of $\alpha$ (see chapter \ref{sub:BifurcationsOfEquilibrium}).
Remarkably, this phenomenon resembles the Eckhaus scenario in spatially
extended diffusive systems \citep{Eckhaus1965,Tuckerman1990}. The
aim of this article is to investigate the transformation of this scenario
under non-local perturbations which destroy the rotational symmetry
of the system. We investigate two different asymptotic cases of small
and large perturbation size $s$. For small $s$, the Eckhaus scenario
persists qualitatively with a modulated Eckhaus line, while for large
$s$, there is a qualitative difference to the unperturbed case that
reflects the new network topology, i.e. the existence of a new loop.

The article is structured as follows. In section~\ref{sec:stability-of-equilibrium}
we discuss the stability of the zero solution (\ref{eq:zero-equilibrium}),
its spectrum and its bifurcations. In section~\ref{sec:Emergent-periodic-orbits}
we find asymptotic expressions for the emerging periodic solutions
in the case of small perturbation size $s$. We also reduce the case
of asymptotically large $s$ to the analysis of an inhomogeneous ring.
Section~\ref{sec:Stability-of-the-psols} deals with the stability
analysis of the periodic solutions and the results are compared to
numerical simulations. Finally, we discuss our findings and give an
outlook on possible extensions and applications of the results in
section \ref{sec:Discussion-1}.

\section{Stability of the zero solution\label{sec:stability-of-equilibrium}}

To study the stability of system (\ref{eq:SL-system}) it is convenient
to identify each variable $z_{j}\left(t\right)\in\mathbb{C}$ with
a two-dimensional real variable $\boldsymbol{z}_{j}\left(t\right)=\left(z_{j,1}\left(t\right),z_{j,2}\left(t\right)\right)^{\top}=\left(\Re\left(z_{j}\left(t\right)\right),\Im\left(z_{j}\left(t\right)\right)\right)^{\top}\in\mathbb{R}^{2}$.
Then, system (\ref{eq:SL-system}) is equivalent to the real system
\begin{eqnarray}
\dot{\boldsymbol{z}}_{j}\left(t\right) & = & \left(M_{\mu}-\left(z_{j,1}^{2}+z_{j,2}^{2}\right)\right)\boldsymbol{z}_{j}+\boldsymbol{z}_{j+1},\nonumber \\
\dot{\boldsymbol{z}}_{N}\left(t\right) & = & \left(M_{\mu}-\left(z_{N,1}^{2}+z_{N,2}^{2}\right)\right)\boldsymbol{z}_{N}+\boldsymbol{z}_{1}+s\boldsymbol{z}_{\ell},\label{eq:SL-system-real}
\end{eqnarray}
where $M_{\mu}=\left[\begin{matrix}\alpha & -\beta\\
\beta & \alpha
\end{matrix}\right]$ is the representation of the multiplication with $\mu$ in $\mathbb{R}^{2\times2}$.

\subsection{Spectrum of the equilibrium}

We linearize system (\ref{eq:SL-system-real}) in $\boldsymbol{z}_{1}=..=\boldsymbol{z}_{N}\equiv0$
and obtain the variational equation 
\begin{align*}
\frac{d}{dt}\delta\boldsymbol{z}_{j}\left(t\right) & =M_{\mu}\delta\boldsymbol{z}_{j}+\delta\boldsymbol{z}_{j+1}\left(t\right),\ j=1,...,N-1,\\
\frac{d}{dt}\delta\boldsymbol{z}_{N}\left(t\right) & =M_{\mu}\delta\boldsymbol{z}_{N}+\delta\boldsymbol{z}_{1}\left(t\right)+s\delta\boldsymbol{z}_{\ell}\left(t\right),
\end{align*}
which can be written as 
\begin{equation}
\frac{d}{dt}\delta Z\left(t\right)=\left[\mathrm{Id}_{N}\otimes M_{\mu}+G_{s}\otimes\mathrm{Id}_{2}\right]\delta Z\left(t\right),\label{eq:variational-eq}
\end{equation}
where $\delta Z=\left(\delta\boldsymbol{z}_{1},...,\delta\boldsymbol{z}_{N}\right)^{\top}$,
$\mathrm{Id}_{N}\in\mathbb{R}^{N\times N}$ is the $N$-dimensional
identity matrix, 
\begin{equation}
G_{s}=\left[\begin{array}{cccc}
0 & 1 &  & 0\\
\vdots & \ddots & \ddots\\
0 &  & \ddots & 1\\
1 & 0 & s & 0
\end{array}\right]\label{eq:coupling-scheme-Gs}
\end{equation}
is the coupling matrix of the perturbed ring, $A\otimes B$ denotes
the tensor product of the two matrices $A$ and $B$. Equation (\ref{eq:variational-eq})
is a simple example of a system, which is treatable by a master stability
function (MSF) approach \citep{Pecora1998}. In our case the MSF $M:\mathbb{C}\longrightarrow\mathbb{R}$
simply reads $M\left(\lambda\right)=\alpha+\Re\left(\lambda\right)$,
where $\lambda$ is an eigenvalue of the coupling matrix $G_{s}$.
Indeed, the spectrum of (\ref{eq:variational-eq}) is 
\begin{equation}
\sigma\left(\mathrm{Id}_{N}\otimes M_{\mu}+G_{s}\otimes\mathrm{Id}_{2}\right)=\sigma\left(M_{\mu}\right)+\sigma\left(G_{s}\right)=\left\{ \mu,\bar{\mu}\right\} +\sigma\left(G_{s}\right).\label{eq:spectrum-equilibrium}
\end{equation}
Taking the real part, we obtain $M\left(\lambda\right)=\alpha+\Re\left(\lambda\right)$.
The spectrum of the coupling matrix $G_{s}$ equals the set of solutions
of the characteristic equation 
\begin{equation}
\chi_{G_{s}}\left(\lambda\right)=\lambda^{N}-s\lambda^{\ell-1}-1=0.\label{eq:CharPoly}
\end{equation}
For $s=0$, $\sigma\left(G_{0}\right)$ consists of the $N$ roots
of unity 
\[
\gamma_{N,k}=e^{i\frac{2\pi k}{N}},\ k=0,...,N-1.
\]
For small $s\ne0$, the roots $\lambda_{k}$, $k=1,...,N$, of eq.~(\ref{eq:CharPoly})
are given (asymptotically) as 
\begin{equation}
\lambda_{k}\left(s\right)=\gamma_{N,k}+\frac{s}{N}\gamma_{N,k}^{\ell}+{\cal O}\left(s^{2}\right),\label{eq:asymp-evs-G}
\end{equation}
as one can readily compute by applying the implicit function theorem
to (\ref{eq:CharPoly}) with base points $(s_{0},\lambda_{0})=(0,\gamma_{N,k})$.
Hence, the spectrum of $G_{s}$ is a weak modulation of the spectrum
of the circulant matrix $G_{0}$ {[}see Fig.~\ref{fig:Spectra}(a),
(b){]}. For large $s$, $\sigma\left(G_{s}\right)$ can be computed
in a similar manner. In Appendix~\ref{sec:Adjacency-spectrum-appendix}
we show that in this case the spectrum of $G_{s}$ splits into two
parts: there are $\ell-1$ roots 
\begin{equation}
\lambda_{1,k}\left(s\right)\approx s^{-\nicefrac{1}{\ell-1}}\gamma_{\ell-1,k},\ k=0,...,\ell-2,\label{eq:ev-inner-circle}
\end{equation}
located close to an inner circle of amplitude $\sim s^{-\nicefrac{1}{\ell-1}}$
and $N-\ell+1$ roots 
\begin{equation}
\lambda_{2,k}\left(s\right)\approx s^{\nicefrac{1}{N-\ell+1}}\gamma_{N-\ell+1,k},\ k=0,...,N-\ell,\label{eq:ev-outer-circle}
\end{equation}
which are close to an outer circle of amplitude $\sim s^{\nicefrac{1}{N-\ell+1}}$
{[}see Fig.~\ref{fig:Spectra}{]}.
\begin{figure}[t]
\centering{}\includegraphics[width=1\textwidth]{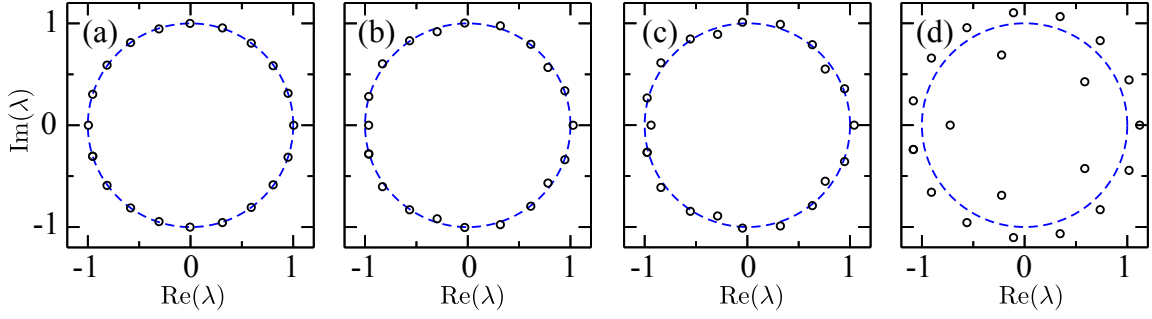}\caption{\label{fig:Spectra}Spectra of the coupling matrix $G_{s}$ (\ref{eq:coupling-scheme-Gs})
for $N=20$ oscillators, a shortcut at node $\ell=6$ and for different
coupling strengths $s$: a) $s=0.1$, b) $s=0.6$, c) $s=1$, d) $s=5$.}
\end{figure}

\subsection{Bifurcations of the equilibrium\label{sub:BifurcationsOfEquilibrium}}

From the formula (\ref{eq:spectrum-equilibrium}) for the spectrum
of the equilibrium $\boldsymbol{z}_{1}=..=\boldsymbol{z}_{N}\equiv0$,
it follows that it is asymptotically stable iff 
\begin{equation}
\Re\left(\lambda\right)<-\alpha,\ \text{for all }\lambda\in\sigma\left(G_{s}\right).\label{eq:stability-equilibrium}
\end{equation}
When the parameter $\alpha$ is increased starting from a value which
satisfies (\ref{eq:stability-equilibrium}), a bifurcation takes place
whenever for some $\lambda\in\sigma\left(G_{s}\right)$: 
\begin{equation}
\alpha=-\Re\left(\lambda\right).\label{eq:bifurcation-points}
\end{equation}
Note that there is always one purely real eigenvalue $\lambda_{1}\left(s\right)=1+\frac{s}{N}+{\cal O}\left(s^{2}\right)$
of $G_{s}$ which has maximal real part among all eigenvalues of $G_{s}$.
Indeed, for any solution $\lambda=\varrho e^{i\varphi}$, $\varrho\ge0$,
of (\ref{eq:CharPoly}), we have 
\[
1=\left|\varrho^{N}e^{iN\varphi}-s\varrho^{\ell-1}e^{i\left(\ell-1\right)\varphi}\right|\ge\varrho^{N}-s\varrho^{\ell-1}=\chi_{G_{s}}\left(\varrho\right)+1.
\]
Hence, $\chi_{G_{s}}\left(\varrho\right)\le0$, and $\lim_{x\to\infty}\chi_{G_{s}}\left(x\right)=+\infty$
implies that there exists a real solution $\varrho_{0}\ge\varrho$
such that $\chi_{G_{s}}\left(\varrho_{0}\right)=0$. Therefore, for
small $s\ge0$ the equilibrium switches stability at 
\[
\alpha_{1}\left(s\right)\approx-\left(1+\frac{s}{N}\right).
\]
Since we assume $\beta\ne0$, this bifurcation is a Hopf bifurcation
and the emerging periodic orbit has frequency $\beta$ at onset. As
for $s=0$, the bifurcation is supercritical for small $s>0$, because
the cubic term of the corresponding normal form depends continuously
on $s$. Therefore, a branch of stable periodic solutions emerges
and exists for $\alpha>\alpha_{1}\left(s\right)$. A further increase
of $\alpha$ leads to a sequence of Hopf bifurcations which give rise
to $N-1$ branches of periodic solutions. As in the case $s=0$, all
bifurcations are supercritical if $s>0$ is small (by continuity).
The same is true if $s>0$ is large, as we show in Appendix~\ref{sec:Supercriticality-of-HB}.
The latter $N-1$ periodic solutions are initially unstable, inheriting
instability from the steady state. The initial frequency $\omega$
of the emerging periodic solution equals the imaginary part of the
crossing eigenvalue, that is $\omega=\beta+\Im\left(\lambda\right)$
for the corresponding $\lambda\in G_{s}$. For the unperturbed ring
($s=0$) the first $\left\lfloor (N-1)/2\right\rfloor $ ($\left\lfloor x\right\rfloor :=\max\left\{ n\in\mathbb{N}:\,n\le x\right\} $)
branches stabilize when they cross the Eckhaus stabilization line
\citep{Yanchuk2008a} 
\begin{equation}
\frac{1}{N}\vert Z\vert^{2}=\frac{3\alpha}{4}+\sqrt{\left(\frac{\alpha}{4}\right)^{2}+\frac{1}{2}},\label{eq:eckhaus-line-seq0}
\end{equation}
where $\left|Z\right|^{2}$ is the (constant) amplitude of a periodic
solution $Z\left(t\right)=(z_{1}\left(t\right),...,z_{N}\left(t\right))^{T}$.
This observation is in striking analogy with the well known Eckhaus
destabilization in diffusive systems \citep{Eckhaus1965,Tuckerman1990}.
Remarkably, it is also found in this unidirectional system which does
not extend to a spatially extended system in a natural manner. In
section \ref{sec:Stability-of-the-psols} we investigate how the added
shortcut changes this scenario.

\section{Emergent periodic orbits\label{sec:Emergent-periodic-orbits}}

\begin{figure}
\centering{}\includegraphics[width=1\textwidth]{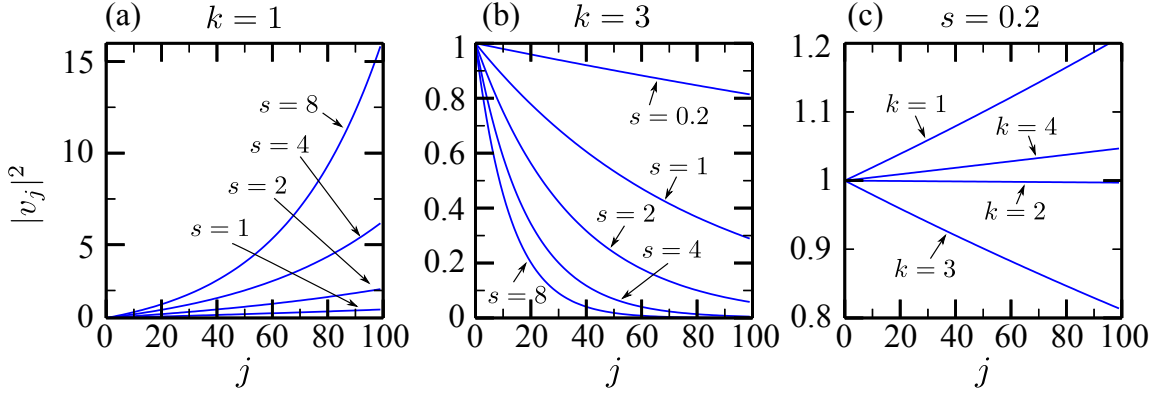}\caption{\label{fig:Profiles}Moduli of the components of the initial profiles
$V\left(0,s\right)=b\left(s\right)$ of emerging periodic solutions
{[}see (\ref{eq:profile-psols}) and (\ref{eq:eigenvector_Gs}){]}
for different wave numbers $k$ (i.e., different eigenvalues $\lambda_{k}\left(s\right)$
of $G_{s}$) and coupling strengths $s$ as indicated in the figure.
For all panels: $N=100$ and $\ell=26$.}
\end{figure}

Let $Z\left(t\right)$ be a periodic solution of (\ref{eq:SL-system}),
which emerges from a Hopf bifurcation at 
\begin{equation}
\alpha\left(s\right)=-\Re\left(\lambda\left(s\right)\right),\label{eq:alpha-k-of-s}
\end{equation}
and which is associated to the eigenvalue $\lambda\left(s\right)$
of the coupling matrix $G_{s}$, see (\ref{eq:asymp-evs-G})--(\ref{eq:ev-outer-circle}).
Because of the $S^{1}$-symmetry of the system, we employ the ansatz
\begin{equation}
Z\left(t\right)=\sqrt{\varepsilon}e^{i\omega\left(\varepsilon,s\right)t}V\left(\varepsilon,s\right),\label{eq:PerAnsatz}
\end{equation}
where 
\begin{equation}
\varepsilon:=\alpha-\alpha\left(s\right)\ge0\label{eq:eps-k-of-s}
\end{equation}
is the parameter distance from the bifurcation point, 
\begin{equation}
V\left(\varepsilon,s\right)=\left(v_{1}\left(\varepsilon,s\right),...,v_{N}\left(\varepsilon,s\right)\right)^{T}\in\mathbb{C}^{N},\label{eq:profile-psols}
\end{equation}
is the profile vector and the frequency $\omega\left(\varepsilon,s\right)$
of $Z\left(t\right)$ is 
\begin{equation}
\omega\left(\varepsilon,s\right)=\beta+\Im\left(\lambda\left(s\right)\right)+{\cal O}\left(\varepsilon\right).\label{eq:freq-psols}
\end{equation}
The emerging periodic orbit is $\varepsilon$-close to the complex
plane spanned by the eigenvector $b\left(s\right)$ of $G_{s}$ corresponding
to $\lambda\left(s\right)$ and tangential at the bifurcation point
itself \citep{Kuznetsov1995}. This means, $V\left(0,s\right)=b\left(s\right)$
with 
\begin{equation}
b\left(s\right)=a\left(s\right)\cdot\left(1,\lambda\left(s\right),\lambda^{2}\left(s\right),...,\lambda^{N-1}\left(s\right)\right)^{T},\label{eq:eigenvector_Gs}
\end{equation}
where one may assume $a\left(s\right)\in\mathbb{R}$ due to the $S^{1}$-symmetry
of the system. Figure~\ref{fig:Profiles} shows several examples
of the profile shapes for different $s$ and $\lambda$. For $\left|\lambda\right|>1$
the emerging solutions become stronger localized at the $N$-th node
$z_{N}$ with increasing $s$ since it scales as $z_{N}\sim\lambda^{N-1}$.
For $\left|\lambda\right|<1$, the localization takes place at $z_{1}\left(t\right)$
for the same reason. 

\subsection{Case I: small perturbation\label{sub:orbits-small-s}}

In this section we consider $s$ to be small. Our aim is to determine
a formal asymptotic expansion for the frequencies (\ref{eq:freq-psols})
and the profiles (\ref{eq:profile-psols}) of the periodic solutions
and to derive evaluable approximate conditions for their stability.
In particular, we are interested in the deformation of the Eckhaus
stabilization line (\ref{eq:eckhaus-line-seq0}). In order to obtain
asymptotic expressions for the profiles in case that $s,\varepsilon>0$,
we linearize the vector field in $\varepsilon=s=0$. For a periodic
solution (\ref{eq:PerAnsatz}), we introduce rescaled, rotating coordinates
$u_{j}\left(t\right)\in\mathbb{C}$, $j=1,...,N$ according to 
\begin{equation}
z_{j}=\sqrt{\varepsilon}e^{i\omega t}v_{j}u_{j},\label{eq:scaled-rotating-coords}
\end{equation}
with $\omega=\omega\left(\varepsilon,s\right)$ and $v_{j}=v_{j}\left(\varepsilon,s\right)$.
Then (\ref{eq:SL-system}) becomes 
\begin{align}
\dot{u}_{j} & =\left(\alpha+i\left(\beta-\omega\right)-\varepsilon\left|v_{j}u_{j}\right|^{2}\right)u_{j}+\frac{v_{j+1}}{v_{j}}u_{j+1},\label{eq:RotatingEquation1}\\
\dot{u}_{N} & =\left(\alpha+i\left(\beta-\omega\right)-\varepsilon\left|v_{N}u_{N}\right|^{2}\right)u_{N}+\frac{v_{1}}{v_{N}}u_{1}+s\frac{v_{\ell}}{v_{N}}u_{\ell},\label{eq:RotatingEquation2}
\end{align}
In rotating coordinates, the equilibrium solution 
\begin{equation}
u_{j}\equiv1,\ j=1,...,N,\label{eq:u-eq-1-sln}
\end{equation}
corresponds to the periodic solution (\ref{eq:PerAnsatz}) of (\ref{eq:SL-system})
and the stability of (\ref{eq:PerAnsatz}) and (\ref{eq:u-eq-1-sln})
is the same. In Appendix \ref{sec:Expansion-of-profiles} we show
that for each eigenvalue $\lambda_{k}\left(s\right)=\gamma_{N,k}+\frac{s}{N}\gamma_{N,k}^{\ell}+{\cal O}\left(s^{2}\right)$
of $G_{s}$, the corresponding branch of periodic solutions has frequencies
\begin{equation}
\omega_{k}\left(\varepsilon,s\right)=\beta+\Im\left(\gamma_{N,k}\right)+\frac{s}{N}\Im\left(\gamma_{N,k}^{\ell}\right)+{\cal O}\left(\left(\left|\varepsilon\right|+\left|s\right|\right)^{2}\right),\label{eq:omega-ebene1}
\end{equation}
and profiles 
\begin{equation}
v_{j}\left(\varepsilon,s\right)=\gamma_{N,k}^{j-1}\left(1+s\frac{j-1}{N}\gamma_{N,k}^{\ell-1}\right)+{\cal O}\left(\left(\left|\varepsilon\right|+\left|s\right|\right)^{2}\right).\label{eq:vjEbene1}
\end{equation}

\subsection{Case II: large perturbation\label{sub:Orbits-large-s}}

In this section we consider $s$ to be large. To treat (\ref{eq:SL-system})
as a weakly perturbed system we perform a change of variables 
\begin{equation}
y_{j}\left(t\right)=\varsigma^{j}z_{j}\left(\varsigma^{2N}t\right),\label{eq:variable-transfo}
\end{equation}
with a small parameter $\varsigma=s^{-\nicefrac{1}{N-\ell+1}}$, which
is the inverse radius of the outer spectral circle of eigenvalues
(\ref{eq:ev-outer-circle}). This transformation normalizes the emerging
profiles (\ref{eq:eigenvector_Gs}) corresponding to the eigenvalues
(\ref{eq:ev-outer-circle}). The transformed variables (\ref{eq:variable-transfo})
satisfy 
\begin{eqnarray}
\dot{y}_{j}\left(t\right) & = & \left(\varsigma^{2N}\mu-\varsigma^{2\left(N-j\right)}\left|y_{j}\left(t\right)\right|^{2}\right)y_{j}\left(t\right)+\varsigma^{2N-1}y_{j}\left(t\right),\ j=1,...,N-1,\nonumber \\
\dot{y}_{N}\left(t\right) & = & \left(\varsigma^{2N}\mu-\left|y_{N}\left(t\right)\right|^{2}\right)y_{N}\left(t\right)+\varsigma^{3N-1}y_{1}+\varsigma^{2N-1}y_{\ell}\left(t\right).\label{eq:transformed-sys-full-large-s}
\end{eqnarray}
Since $\varsigma^{3N-1}y_{1}=\varsigma^{N}\left(\varsigma^{2N-1}y_{1}\right)$
and $\varsigma^{N}$ is small, we consider system (\ref{eq:transformed-sys-full-large-s})
as a small perturbation of 
\begin{eqnarray}
\dot{y}_{j}\left(t\right) & = & \left(\varsigma^{2N}\mu-\varsigma^{2\left(N-j\right)}\left|y_{j}\left(t\right)\right|^{2}\right)y_{j}\left(t\right)+\varsigma^{2N-1}y_{j}\left(t\right),\ j=1,...,N-1,\nonumber \\
\dot{y}_{N}\left(t\right) & = & \left(\varsigma^{2N}\mu-\left|y_{N}\left(t\right)\right|^{2}\right)y_{N}\left(t\right)+\varsigma^{2N-1}y_{\ell}\left(t\right).\label{eq:transformed-sys-large-s-1}
\end{eqnarray}
Although we cannot show that results for the persistence of hyperbolic
invariant manifolds \citep{fenichel1971} apply and assure that (\ref{eq:transformed-sys-large-s-1})
possesses the same hyperbolic invariant manifolds as does (\ref{eq:transformed-sys-full-large-s}),
the truncated system (\ref{eq:transformed-sys-large-s-1}) is a natural
approximation to (\ref{eq:transformed-sys-full-large-s}). Note that
in (\ref{eq:transformed-sys-large-s-1}) the components $y_{1},\,...,\,y_{\ell-1}$
do not couple back to the rest of the system since the weak link from
$y_{1}$ to $y_{N}$ was taken out. Therefore, the dynamics of the
subsystem $y_{\ell},\,...,\,y_{N}$ is independent and, apart from
the zero solution, acts as a periodic force on the attached subsystem
$y_{1},\,...,\,y_{\ell-1}$. In original variables (\ref{eq:transformed-sys-large-s-1})
reads 
\begin{eqnarray}
\dot{z}_{j}\left(t\right) & = & \left(\mu-\left|z_{j}\left(t\right)\right|^{2}\right)z_{j}+z_{j+1}\left(t\right),\ j=1,...,N-1,\nonumber \\
\dot{z}_{N}\left(t\right) & = & \left(\mu-\left|z_{N}\left(t\right)\right|^{2}\right)z_{N}+sz_{\ell}\left(t\right).\label{eq:SL-system-large-s-all-nodes}
\end{eqnarray}
The linearization of (\ref{eq:SL-system-large-s-all-nodes}) at the
zero solution has eigenvalues 
\[
\mu+\nu,\ \,\text{and}\,\ \bar{\mu}+\nu
\]
where $\nu$ is a root of the characteristic equation 
\[
\left(\nu^{N-\ell+1}-s\right)\nu^{\ell-1}=0
\]
of the reduced coupling matrix $H_{s}$, which is obtained by erasing
the link from $z_{1}$ to $z_{N}$ from $G_{s}$. It has eigenvalues
\begin{equation}
\nu=0\ \,\text{and}\ \,\nu=s^{\nicefrac{1}{N-\ell+1}}\gamma_{N-\ell+1,k},\ k=1,...,N-\ell+1.\label{eq:eigenvalue-H-s}
\end{equation}
Periodic orbits which correspond to the algebraically $\left(\ell-1\right)$-fold
eigenvalue $\nu=0$ are localized on the nodes $z_{1},\,...,\,z_{\ell-1}$.

\section{Stability of the periodic orbits\label{sec:Stability-of-the-psols}}

\subsection{Case I: small perturbation\label{sub:stab-small-s}}

\begin{figure}[!t]
\centering{}\centering{}\vspace{-0.6cm}
 \includegraphics[width=1\textwidth]{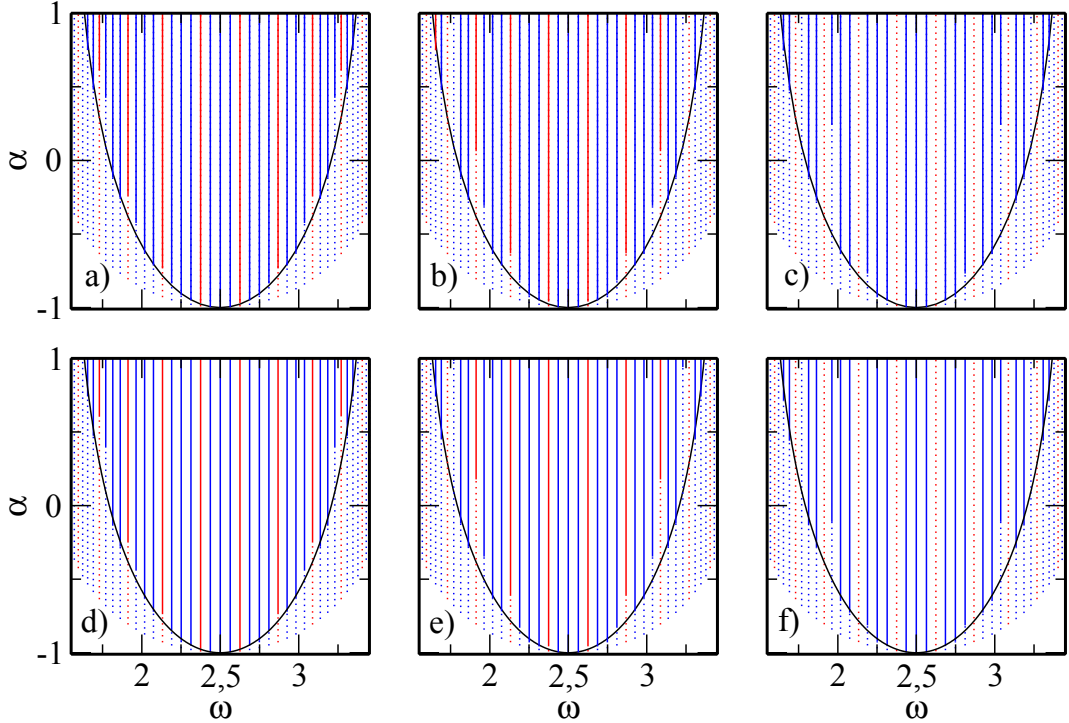}\caption{\label{fig:Eckhaus-freq-vs-a}Bifurcation diagrams for different strengths
$s$ of the shortcut in a ring of $N=100$ oscillators with a shortcut
from node $\ell=26$ to node $N=100$; $\beta=2.5$. Dashed lines
indicate unstable periodic orbits and solid lines stable periodic
orbits. In black the stabilization line for $s=0$ is shown. The frequency
of the solutions is plotted against the bifurcation parameter $\alpha$.
The perturbation strength $s$ increases from left to right: a), d)
$s=0.05$; b), e) $s=0.1$; c), f) $s=0.2$. The upper panels a),
b), c) show the approximated diagram obtained from (\ref{eq:FloquetGlcompact}).
The lower panels d), e), f) show results of numerical bifurcation
analysis of the full system for comparison. These calculations were
carried out with the program AUTO \citep{Doedel2006}.}
\end{figure}

The variational equation for system (\ref{eq:RotatingEquation1})--(\ref{eq:RotatingEquation2})
at the stationary solution (\ref{eq:u-eq-1-sln}) is

\begin{eqnarray*}
\dot{u}_{j} & = & \left(\alpha+\varepsilon+i\left(\beta-\omega\right)-2\varepsilon\left|v_{j}\right|^{2}\right)u_{j}+\frac{v_{j+1}}{v_{j}}u_{j+1},\\
\dot{u}_{N} & = & \left(\alpha+\varepsilon+i\left(\beta-\omega\right)-2\varepsilon\left|v_{N}\right|^{2}\right)u_{N}+\frac{v_{1}}{v_{N}}u_{1}+s\frac{v_{\ell}}{v_{N}}u_{\ell}.
\end{eqnarray*}
We transform the system into $\mathbb{R}^{2N}$ ($\boldsymbol{x}_{j}=(\Re(u_{j}),\Im(u_{j}))^{T}$)
and insert the expansions (\ref{eq:omega-ebene1}) and (\ref{eq:vjEbene1})
to obtain

\begin{eqnarray}
\dot{\boldsymbol{x}}_{j} & = & -\left[M_{\tilde{\lambda}}+2\varepsilon\boldsymbol{\delta}_{11}\right]\boldsymbol{x}_{j}+M_{\tilde{\lambda}}\boldsymbol{x}_{j+1}+{\cal O}\left(\left(\left|\varepsilon\right|+\left|s\right|\right)^{2}\right)\nonumber \\
\dot{\boldsymbol{x}}_{N} & = & -\left[M_{\tilde{\lambda}}+2\varepsilon\boldsymbol{\delta}_{11}\right]\dot{\boldsymbol{x}}_{N}+M_{\tilde{\lambda}}\dot{\boldsymbol{x}}_{1}+sM_{\lambda_{0}^{\ell}}\left[\dot{\boldsymbol{x}}_{\ell}-\dot{\boldsymbol{x}}_{1}\right]+{\cal O}\left(\left(\left|\varepsilon\right|+\left|s\right|\right)^{2}\right)\label{eq:FloquetGl}
\end{eqnarray}
with $\lambda_{0}=\gamma_{N,k}$, $\tilde{\lambda}\left(s\right)=\lambda_{0}+\frac{s}{N}\lambda_{0}^{\ell}$,
the matrix representation $M:\mathbb{C}\to\mathbb{R}^{2\times2}$
{[}as in (\ref{eq:SL-system-real}){]} and $\boldsymbol{\delta}_{mn}=\left(\delta_{jm}\delta_{kn}\right)_{j,k}$
is the matrix which has the entry $1$ at position $\left(j,k\right)$
and zeros everywhere else. We drop the higher order terms in (\ref{eq:FloquetGl})
and write the system in the form 
\begin{equation}
\dot{\boldsymbol{X}}=A\left(\varepsilon,s\right)\boldsymbol{X}=\left[-\mathrm{Id}_{N}\otimes\left(M_{\tilde{\lambda}}+2\varepsilon\boldsymbol{\delta}_{11}\right)+G_{0}\otimes M_{\tilde{\lambda}}+\left(\boldsymbol{\delta}_{N1}-\boldsymbol{\delta}_{N\ell}\right)\otimes sM_{\lambda_{0}^{\ell}}\right]\boldsymbol{X}.\label{eq:FloquetGlcompact}
\end{equation}
Clearly, an MSF approach as in section \ref{sec:stability-of-equilibrium}
is not feasible anymore. However, we have reduced the dynamical problem
to an algebraic one. The eigenvalues of system (\ref{eq:FloquetGlcompact})
can be computed by standard numerical procedures to determine approximately
the stability of the periodic solutions in the vicinity of a bifurcation
point. The eigenvalues of (\ref{eq:FloquetGlcompact}) approximate
the eigenvalues of the exact system (\ref{eq:RotatingEquation1})--(\ref{eq:RotatingEquation2})
at the steady state (\ref{eq:u-eq-1-sln}) up to first order in $\varepsilon$
and $s$. Anyway, this first order approximation leads to good predictions
even for moderate values of the parameters, in particular for $\varepsilon$
{[}see Fig.~\ref{fig:Eckhaus-freq-vs-a}{]}. 

\subsubsection*{Resonances}

An important observation for the perturbed system is that in comparison
to the unperturbed case $s=0$, for small $s>0$ the point of stabilization
of a periodic solution may be altered or not. It remains nearly the
same if the corresponding eigenvalue $\lambda_{k}=\gamma_{N,k}+{\cal O}\left(s\right)$
fulfills 
\begin{equation}
\arg\left(\gamma_{N,k}\right)\approx\arg\left(\gamma_{N,k}^{\ell}\right),\label{eq:phase-matching-cond}
\end{equation}
or equivalently, $\frac{k\left(l-1\right)}{N}\in\mathbb{Z}$. This
corresponds to a situation where both components of the input $z_{1}\left(t\right)+sz_{\ell}\left(t\right)$
to node $z_{N}\left(t\right)$ possess approximately the same phase.
At the same time it is a condition for maximizing the modulus $\left|\lambda_{k}\right|\approx\left|\gamma_{N,k}+\frac{s}{N}\gamma_{N,k}^{\ell}\right|$
and for $\gamma_{N,k}$ to span an eigenmode of both, the unperturbed
cycle of length $N$ and of a cycle which has the same length $N-\ell-1$
as the newly created cycle $\ell\to N\to\left(N-1\right)\to...\to\ell$,
i.e. to be an $N$-th and an $\left(N-\ell+1\right)$-th root of unity.
On the contrary, the antiphase condition 
\begin{equation}
\arg\left(\gamma_{N,k}\right)\approx\arg\left(\gamma_{N,k}^{\ell}\right)+\pi.\label{eq:antiphase-cond}
\end{equation}
causes the point of stabilization to grow rapidly with increasing
$s$, i.e., it destabilizes the corresponding periodic orbit. The
more precisely the equality~(\ref{eq:antiphase-cond}) holds, the
more pronounced is the destabilizing effect of the additional link
{[}see fig.~\ref{fig:Eckhaus-freq-vs-a}{]}.

\subsection{Case II: Large $s$\label{sub:stab-large-s}}

For large $s$ both systems, the original (\ref{eq:SL-system}) and
the truncated (\ref{eq:SL-system-large-s-all-nodes}), admit two types
of periodic solutions emerging in bifurcations corresponding to eigenvalues
of scale $\left|\lambda\right|\approx0$ or $\left|\lambda\right|\approx s^{\nicefrac{1}{N-\ell+1}}$,
respectively {[}see eq.~(\ref{eq:ev-inner-circle}), (\ref{eq:ev-outer-circle}),
and (\ref{eq:eigenvalue-H-s}){]}. In case of system (\ref{eq:SL-system})
all these periodic solutions emerge in form of rotating eigenvectors
(\ref{eq:PerAnsatz})--(\ref{eq:eigenvector_Gs}) of the coupling
matrix $G_{s}$ which leads to a locally pronounced activity in $z_{1}\left(t\right)$
for a corresponding eigenvalue $\left|\lambda\right|\approx0$, and
in $z_{N}\left(t\right)$ for $\left|\lambda\right|\approx s^{\nicefrac{1}{N-\ell+1}}$
{[}see fig.~\ref{fig:Profiles}{]}. The same picture applies for
bifurcations of (\ref{eq:SL-system-large-s-all-nodes}) corresponding
to the simple eigenvalues ($\left|\nu\right|\approx s^{\nicefrac{1}{N-\ell+1}}$)
of $H_{s}$, while the $\ell-1$-fold bifurcation at $\alpha=0$ corresponding
to $\nu=0$ simultaneously creates several solutions which are completely
localized in the attached subsystem $\left(z_{1},...,z_{\ell-1}\right)$
where only the oscillators $z_{1},\,...,\,z_{k}$, $k=1,...,\ell-1$,
are active and all others silent. For each of these solutions the
first active element $z_{k}\left(t\right)=\sqrt{\alpha}e^{i\beta t}$
is located on the limit cycle of an isolated Stuart-Landau oscillator
and all other $z_{j}$, $j<k$, lock either in phase or antiphase
to their input signal $z_{j+1}$. However, these solutions can never
stabilize, since the zero solution of the subsystem $z_{\ell},\,...,\,z_{N}$
is unstable after the first bifurcation corresponding to $\left|\nu\right|\approx s^{\nicefrac{1}{n}}$.
Therefore, it suffices to study the inhomogeneous ring 
\begin{eqnarray}
\dot{z}_{j}\left(t\right) & = & \left(\mu-\left|z_{j}\left(t\right)\right|^{2}\right)z_{j}+z_{j+1}\left(t\right),\ j=1,...,n-1,\nonumber \\
\dot{z}_{n}\left(t\right) & = & \left(\mu-\left|z_{n}\left(t\right)\right|^{2}\right)z_{n}+sz_{1}\left(t\right),\ n=N-\ell+1,\label{eq:SL-system-large-s}
\end{eqnarray}
to understand the possibly stable dynamics of (\ref{eq:SL-system-large-s-all-nodes}).
To approximate the stability of the emerging periodic solutions one
can proceed as for the case of small $s$: write the system in scaled
rotating coordinates (\ref{eq:scaled-rotating-coords}), linearize
around the equilibrium solution (\ref{eq:u-eq-1-sln}), expand the
variational equations in powers of $\varepsilon$ and truncate terms
of order higher than ${\cal O}\left(\varepsilon\right)$.We obtain
the following approximate variational equation {[}see Appendix~\ref{sec:Expansion-of-profiles-inhom}{]}:
\begin{eqnarray}
\dot{\boldsymbol{x}}_{j} & = & -\left[s^{\frac{1}{n}}M_{\gamma_{n,k}}+\varepsilon\left(\left|v_{j}^{0}\left(s\right)\right|^{2}\left(\begin{array}{cc}
3 & 0\\
0 & 1
\end{array}\right)-\left(\begin{array}{cc}
1 & 0\\
0 & 1
\end{array}\right)\right)\right]\boldsymbol{x}_{j}\nonumber \\
 &  & +\left[s^{\frac{1}{n}}M_{\gamma_{n,k}}+\varepsilon\left(\left|v_{j}^{0}\left(s\right)\right|^{2}-1\right)\left(\begin{array}{cc}
1 & 0\\
0 & 1
\end{array}\right)\right]\boldsymbol{x}_{j+1},\label{eq:appr-vareq-large-s-1}
\end{eqnarray}
where $\left|v_{j}^{0}\left(s\right)\right|^{2}=ns^{\frac{2(j-1)}{n}}\frac{s^{\frac{2}{n}}-1}{s^{2}-1}$.
Although the loss of symmetry prevents us from applying an MSF approach,
(\ref{eq:appr-vareq-large-s-1}) enables us to approximate the Floquet
exponents by solving numerically the characteristic equation of (\ref{eq:appr-vareq-large-s-1}).

\begin{figure}[t]
\begin{centering}
\includegraphics[width=0.95\textwidth]{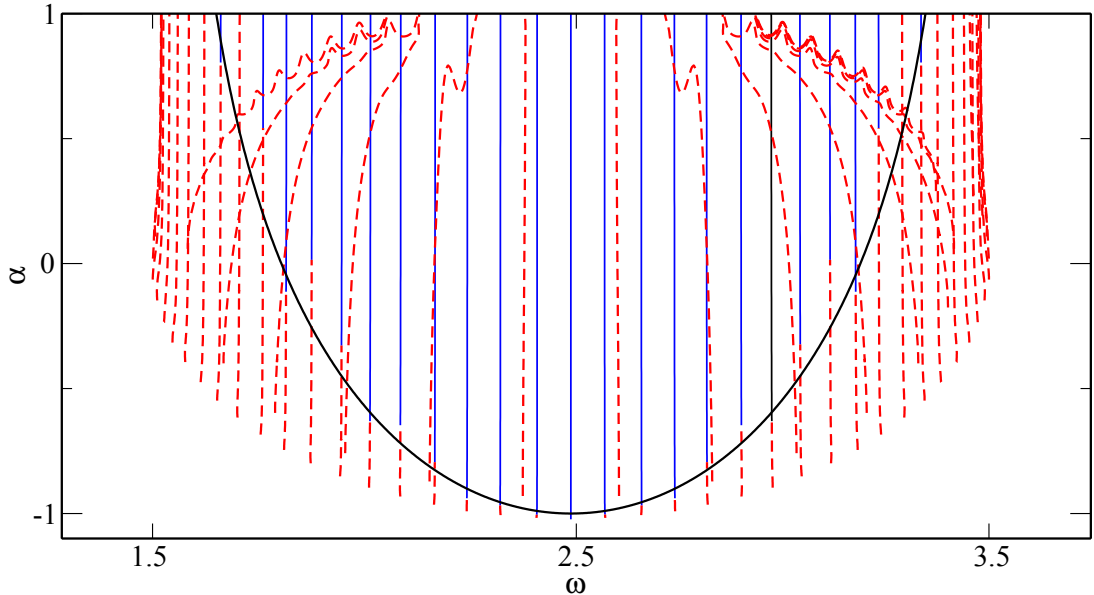}
\par\end{centering}

\caption{\label{fig:bifdiag-freq-vs-a-intermediate-s}Bifurcation diagram for
intermediate strength $s=5$ of the shortcut. Fixed parameters are
$N=100$, $\ell=26$, $\beta=2.5$. Branches of unstable periodic
orbits are indicated by dashed lines and stable branches by solid
lines. Along the black line (\ref{eq:eckhaus-line-seq0}) the stabilization
takes place for $s=0$. The frequency of the solutions is plotted
against the bifurcation parameter $\alpha$. The calculations were
carried out with the program AUTO \citep{Doedel2006}.}
\end{figure}

\section{Discussion\label{sec:Discussion-1}}

We have investigated how the introduction of a shortcut alters the
dynamics in a unidirectional ring of Stuart-Landau oscillators. In
absence of a shortcut, the system exhibits a bifurcation scenario
similar to the Eckhaus instability observed in dissipative media.
For small shortcut strengths $s$ we have found that the Eckhaus stabilization
line is modulated in the following manner: The destabilizing impact
on periodic solutions is stronger for non-resonant modes than for
resonant ones. The latter correspond to wavenumbers that are compatible
with the lengths of both cycles which exist in the perturbed system,
i.e. for the cycle of length $N-\ell+1$ which contains the shortcut,
as well as for the full cycle of length $N$. In contrary to the non-resonant
solutions, the stabilization of the resonant periodic solutions occurs
for similar parameter values as in the case without shortcut. As a
result, one can control the destabilization of a specific set of wavenumbers
via the link position $\ell$ and its strength $s$. In the case of
a large shortcut strength $s$ we have provided an argument that the
cycle of length $N-\ell+1$ dominates the dynamics and stable solutions
can be treated as solutions of a single unidirectional inhomogeneous
ring which has coupling strength $s$ at only one link.

Further investigations will be dedicated to how small perturbations
may be used to select solutions with a specific wavenumber by adding
a corresponding set of resonant links. More generally, our findings
may help to understand how perturbations with a more complicated structure,
consisting of several shortcuts, can influence the dynamics of a unidirectional
ring. Moreover, our observations might even help to locate an unkown
shortcut when one is only allowed to vary a bifurcation parameter
and observe the dynamics, since the modulated Eckhaus line can be
used to identify the shortcut. A strong shortcut can be used in arbitrary
networks in order to localize the activity on the cycles in which
they are contained and amplify the activity in particular on their
targets.

Another point which was not investigated here but deserves a closer
study is the development of profiles far from the bifurcation point.
The simplest, most important phenomenon is that, for increasing $\alpha$,
the exponential tails of the profiles develop into plateaus, where
the profile amplitude is locally constant as a function of the component
index. By taking the limit $N\to\infty$ one can argue that solutions
may consist of several, sharply separated plateaus where each plateau
can possess a different wavenumber. This complies with numerical observations,
although all observed solutions with more than one plateau lie on
the unstable branches which correspond to the inner spectral circle
for large $s$. For sufficiently large values of $s$, those branches
begin to curl with increasing $\alpha$ and seem to be unstable for
arbitrary large values of $\alpha$ {[}see Fig.~\ref{fig:bifdiag-freq-vs-a-intermediate-s}{]}.

\appendix

\section{Adjacency spectrum for large $s$\label{sec:Adjacency-spectrum-appendix}}

To apply the implicit function theorem and continue roots of the characteristic
equation (\ref{eq:CharPoly}) for the case of large $s$, we define
\begin{equation}
F\left(\lambda,\tau,\vartheta\right)=\tau\lambda^{N}-\lambda^{\ell-1}-\vartheta.\label{eq:large-s-function}
\end{equation}
Then $F\left(\lambda,\frac{1}{s},\frac{1}{s}\right)=0$ is equivalent
to (\ref{eq:CharPoly}). We now apply the implicit function theorem
twice to find the two distinct families (\ref{eq:ev-inner-circle})
and (\ref{eq:ev-outer-circle}) of small and large solutions of (\ref{eq:CharPoly}).
Let us first compute 
\begin{eqnarray}
\partial_{\lambda}F\left(\lambda,\tau,\vartheta\right) & = & \tau N\lambda^{N-1}-\left(\ell-1\right)\lambda^{\ell-2},\nonumber \\
\partial_{\tau}F\left(\lambda,\tau,\vartheta\right) & = & \lambda^{N},\qquad\partial_{\vartheta}F\left(\lambda,\tau,\vartheta\right)=-1.\label{eq:large-s-function-derivatives}
\end{eqnarray}
Consider the equation 
\begin{equation}
F\left(\lambda,0,\vartheta\right)=-\lambda^{\ell-1}-\vartheta=0.\label{eq:large-s-small-ev1}
\end{equation}
It possesses $\ell-1$ solutions 
\[
\lambda_{1,k}\left(0,\vartheta\right)=\vartheta^{\frac{1}{\ell-1}}e^{i\frac{\pi}{\ell-1}}\gamma_{\ell-1,k},
\]
$k=0,...,\ell-2$. Assuming (we show that below) that for $\left(\tau,\vartheta\right)\ne0$
one can extend these solutions to smooth functions $\left(\tau,\vartheta\right)\mapsto\lambda_{1,k}\left(\tau,\vartheta\right)$
which solve $F\left(\lambda_{1,k}\left(\tau,\vartheta\right),\tau,\vartheta\right)=0$.
These can be expanded in $\tau=0$ as 
\begin{eqnarray*}
\lambda_{1,k}\left(\tau,\vartheta\right) & = & \lambda_{1,k}\left(0,\vartheta\right)-\partial_{\lambda}F\left(\lambda_{1,k}\left(0,\vartheta\right),0,\vartheta\right)^{-1}\partial_{\tau}F\left(\lambda_{1,k}\left(0,\vartheta\right),0,\vartheta\right)\tau+{\cal O}\left(\tau^{2}\right)\\
 & = & \lambda_{1,k}\left(0,\vartheta\right)+\left(\left(\ell-1\right)\lambda_{1,k}^{\ell-2}\left(0,\vartheta\right)\right)^{-1}\lambda_{1,k}^{N}\left(0,\vartheta\right)\tau+{\cal O}\left(\tau^{2}\right)\\
 & = & \lambda_{1,k}\left(0,\vartheta\right)+\frac{\lambda_{1,k}^{N-\ell+2}\left(0,\vartheta\right)}{\ell-1}\tau+{\cal O}\left(\tau^{2}\right)
\end{eqnarray*}
and, for $\vartheta=\tau$, 
\begin{equation}
\lambda_{1,k}\left(\tau,\tau\right)=\lambda_{1,k}\left(0,\tau\right)+\frac{\lambda_{1,k}^{N-\ell+2}\left(0,\tau\right)}{\ell-1}\tau+{\cal O}\left(\tau^{2}\right).\label{eq:large-s-small-ev-asymp}
\end{equation}
This gives a family of $\ell-1$ solutions of (\ref{eq:CharPoly})
situated near a small circle of radius $\sim\left(1/s\right)^{\frac{1}{\ell-1}}$.
For $\tau\ne0$ , the equation 
\begin{equation}
F\left(\lambda,\tau,0\right)=\tau\lambda^{N}-\lambda^{\ell-1}=0\label{eq:large-s-large-ev1}
\end{equation}
has an $\left(\ell-1\right)$-fold root at $\lambda=0$ which corresponds
to the solutions of (\ref{eq:large-s-small-ev1}) and $N-\ell+1$
roots 
\begin{equation}
\lambda_{2,k}\left(\tau,0\right)=\left|\tau\right|^{-\frac{1}{N-\ell+1}}\gamma_{N-\ell+1,k}\label{eq:Large-s-Large-Eigenvalues}
\end{equation}
$k=0,...,N-\ell$. As before, we obtain an asymptotic representation

\begin{eqnarray}
\lambda_{2,k}\left(\tau,\tau\right) & = & \lambda_{2,k}\left(\tau,0\right)-\partial_{\lambda}F\left(\lambda_{2,k}\left(\tau,0\right),\tau,0\right)^{-1}\partial_{\vartheta}F\left(\lambda_{2,k}\left(\tau,0\right),\tau,0\right)\tau+{\cal O}\left(\tau^{2}\right)\nonumber \\
 & = & \lambda_{2,k}\left(\tau,0\right)+\chi_{G_{s}}^{\prime}\left(\lambda_{2,k}\left(\tau,0\right)\right)^{-1}\tau+{\cal O}\left(\tau^{2}\right)\label{eq:large-s-largel-ev-asymp}
\end{eqnarray}
This is the family of solutions which lie near a larger circle of
radius $\sim s^{\frac{1}{N-\ell+1}}$.

Now let us show that the interval of existence of the implicit functions
$\tau\mapsto\lambda_{1,k}\left(\tau,\vartheta_{0}\right)$ and $\vartheta\mapsto\lambda_{2,k}\left(\tau_{0},\vartheta\right)$,
respectively, indeed contain the points $\tau=\vartheta_{0}$ and
$\vartheta=\tau_{0}$, respectively. Let us denote the implicit function
in question simply $\lambda\left(\tau,\vartheta\right)$. If the implicit
function theorem fails to provide an extension of $\lambda\left(\tau,\vartheta\right)$
in some point $\left(\tau_{\ast},\vartheta_{\ast}\right)>0$, we must
have 
\begin{equation}
\partial_{\lambda}F\left(\lambda_{\ast},\tau_{\ast},\vartheta_{\ast}\right)=\tau_{\ast}N\lambda_{\ast}^{N-1}-\left(\ell-1\right)\lambda_{\ast}^{\ell-2}=0\label{eq:singular-point-ift}
\end{equation}
where $\lambda_{\ast}=\lambda\left(\tau_{\ast},\vartheta_{\ast}\right)$.
Since $F\left(0,\tau,\vartheta\right)=0$ is equivalent to $\vartheta=0$,
we may assume $\lambda_{\ast}\ne0$. Thus, (\ref{eq:singular-point-ift})
is equivalent to 
\begin{equation}
\lambda_{\ast}^{N-\ell+1}=\frac{\ell-1}{N\tau_{\ast}}.\label{eq:ift-ex-cond1}
\end{equation}
Furthermore, from $(\ref{eq:large-s-function})=0$ we obtain $\tau_{\ast}=\lambda^{\ell-1-N}+\vartheta_{\ast}\lambda^{-N}$
which we insert into (\ref{eq:singular-point-ift}) to obtain 
\begin{equation}
\lambda_{\ast}^{\ell-1}=-\frac{N}{N-\ell+1}\vartheta_{\ast}.\label{eq:ift-ex-cond2}
\end{equation}
From (\ref{eq:ift-ex-cond1}) and (\ref{eq:ift-ex-cond2}), we obtain
\[
\tau_{\ast}=\Gamma\left(\vartheta_{\ast}\right):=\frac{\ell-1}{N}\left(\frac{N-\ell+1}{N}\right)^{\frac{N-\ell+1}{\ell-1}}\vartheta_{\ast}^{-\frac{N-\ell+1}{\ell-1}}.
\]
Since $\Gamma:\mathbb{R}_{>0}\to\mathbb{R}_{>0}$ is monotonic, we
have 
\[
\partial_{\lambda}F\left(\lambda\left(\tau,\vartheta_{\ast}\right),\tau,\vartheta_{\ast}\right)\ne0
\]
for all $\tau<\tau_{\ast}=\Gamma\left(\vartheta_{\ast}\right)$. Since
$\Gamma\left(\tau\right)\to\infty$, for $\tau\searrow0$, there exists
$\tau_{0}>0$ such that $\tau<\Gamma\left(\tau\right)$ and $\lambda\left(\tau,\tau\right)$
is defined uniquely, for $0<\tau<\tau_{0}$.

\section{Supercriticality of the Hopf bifurcations\label{sec:Supercriticality-of-HB}}

To show that the bifurcations at $\alpha=-\Re\left(\lambda\right)$
of system (\ref{eq:SL-system-real}) are supercritical for sufficiently
large $s\ge0$, we use the projection method for center manifolds
\citep{Kuznetsov1995}. In the following, we write the vector field
of (\ref{eq:SL-system-real}) as 
\[
f_{\alpha}\left(\boldsymbol{z}\right)=A\boldsymbol{z}+\frac{1}{6}C\left(\boldsymbol{z},\boldsymbol{z},\boldsymbol{z}\right),
\]
where $A=\mathrm{Id}_{N}\otimes M_{\mu}+G_{s}\otimes\mathrm{Id}_{2}$
is the linearization of $f_{\alpha}$ at $\boldsymbol{z}=0$ and the
trilinear function $C$ contains all cubic terms. Furthermore, let
$\boldsymbol{v}\in\mathbb{R}^{2N}$ be an eigenvector of $A$ corresponding
to the eigenvalue $\mu+\lambda$, $\lambda\in\sigma(G_{s})$. (The
case of an eigenvalue $\bar{\mu}+\lambda$ can be treated analogously.)
Let $w\in\mathbb{R}^{2N}$ be the normalized adjoint eigenvector corresponding
to $\boldsymbol{v}$, i.e. $A^{T}\boldsymbol{w}=(\bar{\mu}+\bar{\lambda})\boldsymbol{w}$
and $\left\langle \boldsymbol{w},\boldsymbol{v}\right\rangle =1$.
We have 
\begin{eqnarray}
\boldsymbol{v} & = & \left(1,\lambda,...,\lambda^{N-1}\right)^{T}\otimes\left(i,1\right)^{T},\label{eq:ev}\\
\boldsymbol{w} & = & \frac{1}{\bar{\kappa}}\left(\overline{\lambda}{}^{\ell-1},...,\overline{\lambda},\overline{\lambda}^{N},...,\overline{\lambda}^{\ell}\right)^{T}\otimes\left(i,1\right)^{T},\label{eq:adjoint-ev}
\end{eqnarray}
with $\kappa=2\lambda^{\ell-1}(\left(\ell-1\right)+\left(N-\ell+1\right)\lambda^{N})$.
A Hopf bifurcation at $\alpha$ is supercritical for a negative and
subcritical for a positive first Lyapunov coefficient 
\begin{equation}
l_{1}\left(0\right)=\frac{1}{2\omega_{0}^{2}}\Re\left(\left\langle \boldsymbol{w},C\left(\boldsymbol{v},\boldsymbol{v},\overline{\boldsymbol{v}}\right)\right\rangle \right),\label{eq:first-lyap-coeff}
\end{equation}
where $\omega_{0}=\beta+\Im\left(\lambda\right)\ne0$. Writing $C=(C_{1,1},C_{1,2},\dots,C_{N,1},C_{N,2})$
we obtain 
\begin{equation}
C_{j,1}\left(\boldsymbol{v},\boldsymbol{v},\overline{\boldsymbol{v}}\right)=-6v_{j,1}\left|v_{j,1}\right|^{2}-2\left(2v_{j,1}\left|v_{j,2}\right|^{2}+v_{j,2}^{2}\overline{v}_{j,1}\right)=iC_{j,2}\left(\boldsymbol{v},\boldsymbol{v},\overline{\boldsymbol{v}}\right),\label{eq:components-cubic-form-C}
\end{equation}
for $1\le j\le N$. Using this, we get 
\begin{equation}
\left\langle \boldsymbol{w},C\left(\boldsymbol{v},\boldsymbol{v},\overline{\boldsymbol{v}}\right)\right\rangle =-8\left(\frac{\sum_{j=1}^{\ell-1}\left|\lambda\right|^{2\left(j-1\right)}+\lambda^{N}\sum_{j=\ell}^{N}\left|\lambda\right|^{2\left(j-1\right)}}{\left(\ell-1\right)+\left(N-\ell+1\right)\lambda^{N}}\right)\label{eq:LyapCoeffKompakt}
\end{equation}
\\
 For large $s>0$ we distinguish two cases.For the case $\left|\lambda\right|\sim s^{\nicefrac{1}{N-\ell+1}}$,
we find that 
\[
\Re\left(\left\langle \boldsymbol{w},C\left(\boldsymbol{v},\boldsymbol{v},\overline{\boldsymbol{v}}\right)\right\rangle \right)\to-\infty,\ \text{as }s\to\infty,
\]
and for the case $\left|\lambda\right|\sim s^{-\nicefrac{1}{\ell-1}}$,
\[
\Re\left(\left\langle \boldsymbol{w},C\left(\boldsymbol{v},\boldsymbol{v},\overline{\boldsymbol{v}}\right)\right\rangle \right)\nearrow0,\ \text{as }s\to\infty.
\]
This means that for sufficiently large $s$, we have $l_{1}\left(0\right)<0$
for all eigenvalues. Thus, all bifurcations are supercritical.

\section{Supercriticality for an inhomogeneous ring}

Consider the case $\ell=1$ and $s$ arbitrary, i.e. a ring with inhomogeneous
coupling strengths. Here, using (\ref{eq:LyapCoeffKompakt}) and the
characteristic equation, $\lambda^{N}=1+s$, we get 
\begin{eqnarray}
\left\langle \boldsymbol{w},C\left(\boldsymbol{v},\boldsymbol{v},\overline{\boldsymbol{v}}\right)\right\rangle  & = & -\frac{8\left(1-\left(1+s\right)^{2}\right)}{N\left(1-\left|1+s\right|^{2/N}\right)}.\label{eq:lyap-coeff-inhom-ring}
\end{eqnarray}
Thus, in the case of the inhomogeneous ring the bifurcation is supercritical
for arbitrary $s$.

\section{Expansion of the solution profiles for small perturbations\label{sec:Expansion-of-profiles}}

The linearization of (\ref{eq:RotatingEquation1})--(\ref{eq:RotatingEquation2})
at (\ref{eq:u-eq-1-sln}) is 
\begin{eqnarray}
0 & = & \left(\alpha+i\left(\beta-\omega\right)-\varepsilon\left|v_{j}\right|^{2}\right)+\frac{v_{j+1}}{v_{j}},\label{eq:V-and-omega-full-eqn}\\
0 & = & \left(\alpha+i\left(\beta-\omega\right)-\varepsilon\left|v_{N}\right|^{2}\right)+\frac{v_{1}}{v_{N}}+s\frac{v_{\ell}}{v_{N}},\nonumber 
\end{eqnarray}
From these equations we will now find expressions for the first terms
of the Taylor expansions of the unknown functions 
\[
\omega_{k}\left(\varepsilon,s\right)=\omega_{00}+\varepsilon\omega_{10}+s\omega_{11}+{\cal O}\left(\left(\left|\varepsilon\right|+\left|s\right|\right)^{2}\right)
\]
and 
\[
v_{j}\left(\varepsilon,s\right)=v_{j}^{00}+\varepsilon v_{j}^{10}+sv_{j}^{01}+{\cal O}\left(\left(\left|\varepsilon\right|+\left|s\right|\right)^{2}\right),\ j=1,...,N.
\]

\subsection*{Terms at order ${\cal O}\left(1\right)$}

Considering the terms in $s=\varepsilon=0$ in (\ref{eq:V-and-omega-full-eqn})
yields the circular equations (let $v_{N+1}^{00}:=v_{1}^{00}$) 
\begin{eqnarray}
0 & = & \alpha^{0}+i(\beta-\omega_{00})+\frac{v_{j+1}^{00}}{v_{j}^{00}},\label{eq:circ-eq-O1}
\end{eqnarray}
with the shorthand $\alpha^{0}=\alpha_{k}\left(0\right)=-\cos\left(2\pi k/N\right)$.
This leads to 
\[
\left(-\alpha^{0}-i\left(\beta-\omega_{00}\right)\right)^{N}=1
\]
which contains no new information, since it only determines $\alpha^{0}=-\Re\left(\lambda_{0}\right)$
and $\omega_{00}=\beta+\Im\left(\lambda_{0}\right)$ with an $N$-th
root of unity 
\[
\lambda_{0}=\lambda_{k}\left(0\right)=e^{i\frac{2\pi k}{N}}.
\]
For the profile, (\ref{eq:circ-eq-O1}) yields 
\begin{equation}
v_{j+1}^{00}=\lambda_{0}^{j}r_{0},\label{eq:profile-O1}
\end{equation}
with a hitherto unknown scale factor $r_{0}:=v_{1}^{00}$. Without
loss of generality, one may choose $v_{1}\left(\varepsilon,s\right)\in\mathbb{R}_{\ge0}$,
because of the phase shift invariance of (\ref{eq:PerAnsatz}). In
particular, we then have $r_{0}\in\mathbb{R}_{+}$.

\subsection*{Terms at order ${\cal O}\left(\varepsilon\right)$}

At first order in $\varepsilon$ we obtain another set of circular
equations (let again $v_{N+1}^{10}=v_{1}^{10}$) 
\begin{eqnarray*}
0 & = & 1-i\omega_{10}-r_{0}^{2}+\frac{\lambda_{0}}{r_{0}}\left(\frac{v_{j+1}^{10}}{\lambda_{0}^{j}}-\frac{v_{j}^{10}}{\lambda_{0}^{j-1}}\right)
\end{eqnarray*}
which leads us to the following recursion 
\begin{equation}
\frac{v_{j+1}^{10}}{\lambda_{0}^{j}}=\left(r_{0}^{2}-\left(1-i\omega_{10}\right)\right)\frac{r_{0}}{\lambda_{0}}+\frac{v_{j}^{10}}{\lambda_{0}^{j-1}}\label{eq:Rekursion1}
\end{equation}
Defining $a_{j}=\frac{v_{j}^{10}}{\lambda_{0}^{j-1}}$ and $A=\left(r_{0}^{2}-(1-i\omega_{10})\right)\frac{r_{0}}{\lambda_{0}}$,
equation (\ref{eq:Rekursion1}) can be written as $a_{j+1}=A+a_{j}$
with the solution $a_{j+1}=jA+a_{1}$. Because of the circularity
$a_{N+1}=a_{1}$ we then have $a_{1}=NA+a_{1}$ which determines $A=0$.
Therefore, $r_{0}^{2}=(1-i\omega_{10})$ and finally 
\begin{equation}
\omega_{10}=0\ \text{and}\ r_{0}=1\label{eq:omega10-and-r0}
\end{equation}
That means at first order the frequency of the oscillations does not
depend on $\varepsilon$. For the profiles $v_{j}^{10}$ it follows
(as at order ${\cal O}\left(1\right)$) 
\begin{equation}
v_{j+1}^{10}=\lambda_{0}^{j}r_{1\text{0}},\label{eq:profile-O-eps}
\end{equation}
with $r_{10}:=v_{1}^{10}\in\mathbb{R}$.

\subsection*{Terms at $O\left(s\right)$}

Due to the fact that $s$ perturbs the network symmetry, the equations
at order ${\cal O}\left(s\right)$ are more complex than those at
order ${\cal O}\left(\varepsilon\right)$. Again, the linear terms
in $s$ of (\ref{eq:V-and-omega-full-eqn}) give a recursive formula
\begin{eqnarray}
\frac{v_{j+1}^{01}}{\lambda_{0}^{j}} & = & \frac{1}{\lambda_{0}}\left(i\omega_{01}-\alpha^{1}\right)+\frac{v_{j}^{01}}{\lambda_{0}^{j-1}},\nonumber \\
v_{1}^{01} & = & \frac{1}{\lambda_{0}}\left(i\omega_{01}-\alpha^{1}\right)-\lambda_{0}^{\ell-1}+\frac{v_{N}^{01}}{\lambda_{0}^{N-1}},\label{eq:v01Relation}
\end{eqnarray}
where $\alpha^{1}=-\frac{d}{ds}\Re\left(\lambda_{k}\left(s\right)\right)\mid_{s=0}$.
For $j<N$ we have $\frac{v_{j+1}^{01}}{\lambda_{0}^{j}}=\frac{j}{\lambda_{0}}\left(i\omega_{01}-\alpha^{1}\right)+v_{1}^{01}$.
Inserting the resulting expression for $v_{N}^{01}$ in the second
equation gives $v_{1}^{01}=\frac{N}{\lambda_{0}}\left(i\omega_{01}-\alpha^{1}\right)-\lambda_{0}^{\ell-1}+v_{1}^{01}$.
Therefore, $N\left(i\omega_{01}-\alpha^{1}\right)=\lambda_{0}^{\ell}$
or equivalently 
\begin{equation}
\omega_{01}=\frac{1}{N}\Im\left(\lambda_{0}^{\ell}\right)\ \text{and}\ \alpha^{1}=-\frac{1}{N}\Re\left(\lambda_{0}^{\ell}\right).\label{eq:omega-and-alpha-O-s}
\end{equation}
The perturbations $v_{j+1}^{01}$ of the profile are then determined
as 
\begin{eqnarray}
v_{j+1}^{01} & = & \left(\frac{j}{N}\lambda_{0}^{\ell-1}+r_{01}\right)\lambda_{0}^{j}.\label{eq:profile-O-s}
\end{eqnarray}
up to a scaling $r_{01}:=v_{1}^{01}\in\mathbb{R}$ as above.

\subsection*{Higher order terms}

To determine the amplitudes $r_{10}$ and $r_{01}$ of the perturbations
$v_{j}^{10}$ and $v_{j}^{01}$, $j=1,...,N$, we need to calculate
the second order terms ${\cal O}\left(\varepsilon^{2}\right)$, ${\cal O}\left(s^{2}\right)$,
and ${\cal O}\left(\varepsilon\cdot s\right)$ of (\ref{eq:V-and-omega-full-eqn})
due to the nonlinear term. We omit this here and just give the resulting
values 
\begin{equation}
r_{10}=0\ \text{and}\ r_{01}=0.\label{eq:r10-and-r01}
\end{equation}
The vanishing $r_{10}$ means that at first order the profiles do
not depend on $\varepsilon$.

\section{Expansion of the solution profiles for the inhomogeneous ring \label{sec:Expansion-of-profiles-inhom} }

A periodic solution of (\ref{eq:SL-system-large-s}) corresponds to
a fixed point in co-rotating coordinates (\ref{eq:scaled-rotating-coords}).
The stability of this fixed point is governed by its variational equations
\begin{eqnarray}
\dot{\boldsymbol{x}}_{j} & = & \left(M_{\alpha+i\left(\beta-\omega\right)}-\varepsilon\left(s\right)\left|v_{j}\right|^{2}\left[\begin{matrix}3 & 0\\
0 & 1
\end{matrix}\right]\right)\dot{\boldsymbol{x}}_{j}+M_{\frac{v_{j+1}}{v_{j}}}\dot{\boldsymbol{x}}_{j+1},\label{eq:vareqn-psol-large-s-11}\\
\dot{\boldsymbol{x}}_{n} & = & \left(M_{\alpha+i\left(\beta-\omega\right)}-\varepsilon\left(s\right)\left|v_{n}\right|^{2}\left[\begin{matrix}3 & 0\\
0 & 1
\end{matrix}\right]\right)\dot{\boldsymbol{x}}_{n}+sM_{\frac{v_{1}}{v_{n}}}\dot{\boldsymbol{x}}_{1},\label{eq:vareqn-psol-large-s-12}
\end{eqnarray}
Higher order terms of $\omega=\omega\left(\varepsilon,s\right)$ and
$v_{j}=v_{j}\left(\varepsilon,s\right)$ can be determined by the
equations 
\begin{eqnarray}
0 & = & \left(\alpha+i\left(\beta-\omega\right)\right)-\varepsilon\left|v_{j}\right|^{2}+\frac{v_{j+1}}{v_{j}},\label{eq:V-and-omega-inh-1}\\
0 & = & \left(\alpha+i\left(\beta-\omega\right)\right)-\varepsilon\left|v_{n}\right|^{2}+s\frac{v_{1}}{v_{n}},\label{eq:V-and-omega-inh-2}
\end{eqnarray}
which are obtained from inserting the solution ansatz (\ref{eq:PerAnsatz})
into (\ref{eq:SL-system-large-s}). Solving for real and imaginary
parts yields the conditions 
\begin{eqnarray}
\omega-\beta & = & \Im\left(\frac{v_{j+1}}{v_{j}}\right)=s\Im\left(\frac{v_{1}}{v_{n}}\right),\label{eq:cond-frequency}\\
\alpha & = & \varepsilon\left|v_{j}\right|^{2}-\Re\left(\frac{v_{j+1}}{v_{j}}\right)=\varepsilon\left|v_{n}\right|^{2}-s\Re\left(\frac{v_{1}}{v_{n}}\right).\label{eq:cond-profile}
\end{eqnarray}
Note that we expand the unknown functions only in $\varepsilon$,
keeping $s$ arbitrary. Using (\ref{eq:cond-frequency}) and (\ref{eq:cond-profile})
the variational equations (\ref{eq:vareqn-psol-large-s-11})--(\ref{eq:vareqn-psol-large-s-12})
write 
\begin{eqnarray}
\dot{\boldsymbol{x}}_{j} & = & -\left[M_{\frac{v_{j+1}}{v_{j}}}+2\varepsilon\mid v_{j}\mid^{2}\left(\begin{array}{cc}
1 & 0\\
0 & 0
\end{array}\right)\right]\boldsymbol{x}_{j}+M_{\frac{v_{j+1}}{v_{j}}}\boldsymbol{x}_{j+1},\label{eq:vareqn-psol-large-s-21}\\
\dot{\boldsymbol{x}}_{n} & = & -\left[sM_{\frac{v_{1}}{v_{n}}}+2\varepsilon\mid v_{n}\mid^{2}\left(\begin{array}{cc}
1 & 0\\
0 & 0
\end{array}\right)\right]\boldsymbol{x}_{n}+sM_{\frac{v_{1}}{v_{n}}}\boldsymbol{x}_{1}.\label{eq:vareqn-psol-large-s-22}
\end{eqnarray}
In the rest of this section, we find expressions for $\frac{v_{j+1}}{v_{j}}\left(0,s\right)$
and $\frac{\partial}{\partial\varepsilon}\left[\frac{v_{j+1}}{v_{j}}\left(\varepsilon,s\right)\right]_{|\varepsilon=0}$
which can be inserted into (\ref{eq:vareqn-psol-large-s-21})--(\ref{eq:vareqn-psol-large-s-22})
to obtain (\ref{eq:appr-vareq-large-s-1}). Let 
\[
\omega\left(\varepsilon,s\right)=\omega_{0}\left(s\right)+\varepsilon\omega_{1}\left(s\right)+{\cal O}\left(\varepsilon^{2}\right)
\]
and 
\[
v_{j}\left(\varepsilon,s\right)=v_{j}^{0}\left(s\right)+\varepsilon v_{j}^{1}\left(s\right)+{\cal O}\left(\varepsilon^{2}\right),\ j=1,...,N.
\]
For $\varepsilon=0$, equations (\ref{eq:V-and-omega-inh-1})--(\ref{eq:V-and-omega-inh-2})
yield 
\begin{eqnarray*}
v_{j+1}^{0} & = & \left(i\left(\omega_{0}-\beta\right)-\alpha_{0}\right)v_{j}^{0}=...=\left(i(\omega_{0}-\beta)-\alpha_{0}\right)^{j}v_{1}^{0},\\
v_{1}^{0} & = & \frac{1}{s}\left(i\left(\omega_{0}-\beta\right)-\alpha_{0}\right)v_{n}^{0},
\end{eqnarray*}
where $\alpha_{0}=-\Re\left(\lambda_{k}\left(s\right)\right)=-s^{\nicefrac{1}{n}}\cos\left(2\pi k/n\right)$
is the critical value at which the periodic solution emerges, $\omega_{0}=\beta+\Im\left(\lambda_{k}\right)$
is the initial frequency, and the initial profile is given by $v_{j}^{0}=\lambda_{k}^{j-1}v_{1}^{0}$.

For terms of order ${\cal O}\left(\varepsilon\right)$ one obtains
\[
\omega_{1}=0,\ \left|v_{1}^{0}\right|^{2}=n\frac{s^{\frac{2}{n}}-1}{s^{2}-1},\ \text{and }v_{j}^{1}=v_{j}^{0}\sum_{l=0}^{j-1}\left(\vert v_{j-l}^{0}\vert^{2}-1\right)+\lambda_{k}^{j}v_{1}^{1}
\]
To approximate (\ref{eq:vareqn-psol-large-s-21})--(\ref{eq:vareqn-psol-large-s-22}),
we calculate the first order terms of the quotients $v_{j+1}/v_{j}$,
which gives at ${\cal O}(\varepsilon^{0})$ 
\[
\frac{v_{j+1}^{0}}{v_{j}^{0}}=s^{\nicefrac{1}{n}}\gamma_{n,k},
\]
and at ${\cal O}(\varepsilon)$ 
\begin{eqnarray*}
\left(\frac{v_{j+1}}{v_{j}}\right)^{1} & = & \left|v_{j}^{0}\right|^{2}-1.
\end{eqnarray*}

\bibliographystyle{unsrt}
\addcontentsline{toc}{section}{\refname}\bibliography{DISS}

\end{document}